\documentclass[11pt, a4paper]{amsart}

\usepackage[T1]{fontenc}
\usepackage{ae}

\renewcommand\AA{\mathbb{A}}
\newcommand\OO{\mathcal{O}}
\newcommand\QQ{\mathbb{Q}}
\newcommand\ZZ{\mathbb{Z}}
\newcommand\FF{\mathbb{F}}
\newcommand\PP{\mathbb{P}}
\newcommand\LL{\mathcal{L}}
\newcommand\DD{\mathcal{D}}
\newcommand\KK{\mathbb{K}}
\newcommand\Ptwo{{\PP^2}}
\newcommand\Pthree{{\PP^3}}
\newcommand\inj{\hookrightarrow}
\newcommand\rto{\dashrightarrow}
\newcommand{\Aone}{{\mathbf A}_1}
\newcommand{\Atwo}{{\mathbf A}_2}

\newcommand{\Afour}{{\mathbf A}_4}
\newcommand{\Esix}{{\mathbf E}_6}
\newcommand{\Eseven}{{\mathbf E}_7}
\newcommand{\Eeight}{{\mathbf E}_8}
\newcommand{\Dfive}{{\mathbf D}_5}
\newcommand{\anti}[1]{-K_{S_#1}}
\newcommand{\tanti}[1]{2\cdot(\anti #1)}
\DeclareMathOperator{\Pic}{Pic}
\DeclareMathOperator{\Cox}{Cox}
\DeclareMathOperator{\Spec}{Spec}
\DeclareMathOperator{\rad}{rad}
\newcommand{\qa}[1]{q_{#1}}
\newcommand{\qb}[1]{q'_{#1}}
\newcommand{\qc}[1]{q''_{#1}}

\newtheorem{theorem}{Theorem}
\newtheorem{lemma}[theorem]{Lemma}
\theoremstyle{definition}
\newtheorem{definition}[theorem]{Definition} 
\theoremstyle{remark}
\newtheorem{remark}[theorem]{Remark}
\newtheorem{example}[theorem]{Example}

\begin{document}

\title{On the Cox ring of Del Pezzo surfaces}
\author{Ulrich Derenthal} 
\address{Mathematisches Institut, Universit\"at G\"ottingen,
  Bunsenstr. 3-5, 37073 G\"ottingen, Germany}
\email{derentha@math.uni-goettingen.de} 
\date{March 4, 2006}

\begin{abstract}
  Let $S_r$ be the blow-up of $\Ptwo$ in $r$ general points, i.e., a
  smooth Del Pezzo surface of degree $9-r$. For $r \le 7$, we
  determine the quadratic equations defining its Cox ring explicitly.
  The ideal of the relations in $\Cox(S_8)$ is calculated up to
  radical. As conjectured by Batyrev and Popov, all the generating
  relations are quadratic.
\end{abstract}

\maketitle

\tableofcontents

\section{Introduction}

Over an algebraically closed field $\KK$, the blow-up $S_r$ of $\Ptwo$ in $r$
points in general position\footnote{I.e., no three points on one line, no six
  points on a conic, no eight points with one of them a double point on a
  cubic.} ($1 \le r \le 8$) is a Del Pezzo surface of degree $d = 9-r$. Its
Picard group is a free $\ZZ$-module of rank $r+1$.

Once we have chosen representatives $\LL_0, \dots, \LL_r$ for a basis of
$\Pic(S_r)$, we can define its \emph{Cox ring}, or \emph{total coordinate
  ring}, as
\[\Cox(S_r) := \bigoplus_{(\nu_0, \dots, \nu_r) \in \ZZ^{r+1}}
\Gamma(S_r,\LL_0^{\otimes \nu_0} \otimes \dots \otimes \LL_r^{\otimes
  \nu_r}).\] The multiplication of sections induces the multiplication in
$\Cox(S_r)$.  The Cox ring is graded by $\Pic(S_r)$ and is independent of the
choice of the basis.

The \emph{intersection form} is a non-degenerate bilinear form on $\Pic(S_r)$.
We will write it as $(D_1,D_2)$ for $D_1,D_2 \in \Pic(S_r)$. (We will often
use the same notation for divisors and their class in $\Pic(S_r)$. It will be
clear from the context what is meant.) A prime divisor $D$ whose
self-intersection number $(D,D)$ is negative is called a \emph{negative
  curve}. On smooth Del Pezzo surfaces, every negative curve has
self-intersection number $-1$.

For $r \in \{3, \dots, 7\}$, $\Cox(S_r)$ is generated by non-zero sections of
the $N_r$ negative curves (\cite[Theorem 3.2]{MR2029863}), see
Table~\ref{tab:relations} for the values of $N_r$. For $r=8$, we must add two
independent sections of $\Gamma(S_8, \anti 8)$. Let $R_r$ be the free
polynomial ring whose variables correspond to these generators of $\Cox(S_r)$.
We want to determine the relations between these generators.

For $r \le 3$, the Cox ring is a polynomial ring in $r+3$ generators. This is
due to the fact that in these cases, $S_r$ is toric (see \cite{MR95i:14046}
for Cox rings of toric varieties).

\begin{definition}
  For $n \ge 1$, a divisor class $D$ is called an \emph{$(n)$-ruling}
  if $D=D_1+D_2$ for two negative curves $D_1,D_2$ whose intersection
  number $(D_1,D_2)$ is $n \ge 1$. A $(1)$-ruling is also called a
  \emph{ruling}.
\end{definition}

Each $(n)$-ruling defines quadratic relations between 
generators of $\Cox(S_r)$, see Lemma~\ref{lem:quadric_relations}. 
Relations coming from $(1)$-rulings define an ideal $I_r \subset R_r$.
For $r \in \{4,5,6\}$, $\Cox(S_r) = R_r / \rad(I_r)$ by \cite[Theorem
4.9]{MR2029863}. We extend this result to $r\in \{7,8\}$
as follows:

\begin{theorem}\label{thm:relations}
  For $r \in \{4, \dots, 8\}$, we have $\Cox(S_r) = R_r / \rad(J_r)$, where
  \begin{itemize}
  \item for $r \in \{4, 5, 6\}$, $J_r:=I_r$;
  \item the ideal $J_7$ is generated by the 504 quadratic relations
    coming from the 126 rulings, and 25 quadratic relations coming
    from the $(2)$-ruling $\anti 7$;
  \item the ideal $J_8$ is generated by the 10800 quadratic relations
    coming from the 2160 rulings, 6480 quadratic
    relations coming from 240 $(2)$-rulings, and 119 quadratic
    relations coming from the $(3)$-ruling $\tanti 8$.
  \end{itemize}
\end{theorem}

It is known that the ideal $I_4$ is radical (see \cite{MR2029863}). 
Batyrev proved that the same holds for $I_5$ (unpublished). Here we prove:

\begin{theorem}\label{thm:radical}
  For $r \in \{4, \dots, 7\}$, the ideals $J_r$ are radical, and \[\Cox(S_r) =
  R_r / J_r.\]
\end{theorem}

It was conjectured by Batyrev and Popov that the ideal of relations
defining $\Cox(S_r)$ is generated by quadrics for $r \in \{4, \dots,
8\}$, see \cite[Conjecture 4.3]{MR2029863}. To prove this conjecture, 
it now remains to show that $J_8$ is radical.

After recalling some general results on Del Pezzo surfaces in Section
\ref{sec:generalities}, we will handle the cases $r \in \{6, 7,
8\}$ separately. 

\

\noindent{\textbf{Acknowledgments.}} I am grateful to V. Batyrev for providing 
me with similar calculations for the case of degree 4 Del Pezzo
surfaces. I thank H.-C. Graf v. Bothmer for help with the calculation
of the quadratic relations in case of the cubic surface.

\section{Smooth Del Pezzo surfaces}\label{sec:generalities}

In this section, we summarize some facts on smooth Del Pezzo surfaces.

\begin{itemize}
\item Let $E_1, \dots, E_r$ be the exceptional divisors of the blow-up of
  $\mathbb P^2$ in $r$ points $p_1, \dots, p_r$ in general position. A basis
  of $\Pic(S_r)$ is given by (the classes of) $H, E_1, \dots, E_r$, where $H$
  is the pullback of the hyperplane section in $\mathbb P^2$.
\item In terms of this basis, the intersection form is given by a diagonal
  matrix of size $r+1$ whose diagonal is $(1, -1, \dots, -1)$. The
  anticanonical divisor is $\anti r = 3H-(E_1+\dots+E_r)$.
\item The curves with self-intersection number $-1$ are described in
  \cite[Theorem 2.1]{MR2029863}. There are no curves whose
  self-intersection is $\le -2$.
\item The Weyl group $W_r$ acting on $\Pic(S_r)$ depends on $r$:
  \[\begin{array}[h]{|c||c|c|c|c|c|c|c|c|}
    \hline
    r & 1 & 2 & 3 & 4 & 5 & 6 & 7 & 8 \\
    \hline
    W_r & \Aone & \Atwo & \Atwo+\Aone & \Afour & \Dfive & \Esix & 
    \Eseven & \Eeight\\
    \hline
  \end{array}\]
  For more details, see \cite[Section 2]{MR2029863}.
\end{itemize}

As explained in the introduction, for $r \le 6$, all relations in the
Cox ring are induced by rulings, and these relations also play an
important role for $r \in \{7,8\}$. More precisely, by the discussion
following \cite[Remark 4.7]{MR2029863},
each ruling is represented in $r-1$ different ways as the sum of two
negative curves, giving $r-3$ linearly independent quadratic relation
in $\Cox(S_r)$. Therefore, if each of the $N_r$ negative curves
intersects $n_r$ negative curves with intersection number $1$, we have
$N'_r = (N_r \cdot n_r)/2$ pairs, the number of rulings is $N''_r =
N'_r/(r-1)$, and the number of quadratic relations coming from rulings
is $N''_r\cdot (r-3)$ (see Table~\ref{tab:relations}).
\begin{table}[h]
  \centering
  \[\begin{array}[h]{|c||c|c|c|c|c|c|c|c|c|c|}
    \hline
    r & 3 & 4 & 5 & 6 & 7 & 8 \\
    \hline\hline
    N_r & 6 & 10 & 16 & 27 & 56 & 240\\
    n_r & 2 & 3 & 5 & 10 & 27 & 126\\
    N''_r & 3 & 5 & 10 & 27 & 126 & 2160\\
    \hline\hline
    \text{relations} & 0 & 5 & 20 & 81 & 504 & 10800\\
    \hline
  \end{array}\] 
  \smallskip
  \caption{The number of relations coming from rulings.}
  \label{tab:relations}
\end{table}

Now we describe how to obtain explicit equations for $\Cox(S_r)$ and
how to prove Theorem~\ref{thm:relations} and Theorem~\ref{thm:radical}.
We isolate the steps that must be carried out for each of the degrees
3, 2, and 1 and complete the proofs in the following sections.

\textbf{Choice of coordinates.}  Choose coordinates for 
$p_1, \dots, p_r \in \Ptwo$. We may assume that the first four points are
\begin{equation}\label{eq:four_points}
p_1 = (1:0:0), \quad p_2 = (0:1:0), \quad p_3 = (0:0:1), \quad p_4 = (1:1:1).
\end{equation}
By the general position requirement, the other points must have non-zero
coordinates, and we can write $p_j = (1:\alpha_j:\beta_j)$ for $j \in \{5,
\dots, r\}$.

\textbf{Curves in $\Ptwo$.} As explained in the introduction, $\Cox(S_r)$ is
generated by sections of the negative curves for $r \le 7$. For a negative
curve $D$, we denote the corresponding section by $\xi(D)$, and for a
generating section $\xi$, let $D(\xi)$ be the corresponding divisor. For $r =
8$, we need two further generators: linearly independent sections $\kappa_1,
\kappa_2$ of $\anti 8$. Let $K_1 := D(\kappa_1)$, $K_2 := D(\kappa_2)$ be the
corresponding divisors in the divisor class $\anti 8$.

Let $\DD_r$ be the set of divisors corresponding to sections
generating $\Cox(S_r)$ (including $K_1, K_2$ if $r=8$).

We need an explicit description of the image of each generator $D$ of
$\Cox(S_r)$ under the projection $\pi: S_r \to \Ptwo$.  According to the seven
cases in \cite[Theorem 2.1]{MR2029863}, $\pi(D)$ can be a curve, determined by
a form $f_D$ of degree $d \in \{1, \dots, 6\}$, or a point (if $D = E_i$). If
$\pi(D)$ is a point, the convention to choose $f_D$ as a non-zero constant
will be useful later.

For $r = 8$, we have the following situation: The image of $K_i$ is a cubic
through the eight points $p_1, \dots, p_8$. The choice of two linearly
independent sections $\kappa_1, \kappa_2$ corresponds to the choice of two
independent cubic forms $f_{K_1}, f_{K_2}$ vanishing in the eight points.
Every cubic through these points has the form $a_1f_{K_1}+a_2f_{K_2}$ where
$(a_1,a_2) \ne (0,0)$, and the cubic does not change if we replace $(a_1,a_2)$
be a non-zero multiple. This gives a one-dimensional projective space of
cubics through the eight points.

Let $X_1, \dots, X_n$ be the monomials of degree $d$ in three variables
$x_0,x_1,x_2$. For $D \in \DD_r$, we can write
\[f_D = \sum_{i=1}^n a_i\cdot X_i\] for suitable coefficients $a_i$, which we
can calculate in the following way: If $p_j$ lies on $\pi(D)$, this gives a
linear condition on the coefficients $a_i$ by substituting the coordinates of
$p_j$ for $x_0,x_1,x_2$. If $p_j$ is a double point of $\pi(D)$, all partial
derivatives of $f_D$ must vanish at this point, giving three more linear
conditions. If $p_j$ is a triple point, we get six more linear conditions from
the second derivatives. With $p_1, \dots, p_r$ in general position, we check
that these conditions determine $f_D$ uniquely up to a non-zero constant.

\textbf{Relations corresponding to $(n)$-rulings.} Suppose that an
$(n)$-ruling $D$ can be written as $D_j+D'_j$ for $k$ different pairs $D_j,
D'_j \in \DD_r$ where $j \in \{1, \dots, k\}$. Then the products \[f_{D_1}\cdot
f_{D'_1}, \dots, f_{D_k}\cdot f_{D'_k}\] are $k$ homogeneous forms of the same
degree $d$, and they span a vector space of dimension $n+1$ in the space of
homogeneous polynomials of degree $d$. Therefore, there are $k-(n+1)$
independent relations between them, which we write as
\[\sum_{j=1}^k a_{j,i} \cdot f_{D_j}\cdot f_{D'_j}=0 \qquad\text{for
  $i \in \{1, \dots, k-(n+1)\}.$}\] for suitable constants $a_{j,i}$. They
give an explicit description of the quadric relations coming from $D$:

\begin{lemma}\label{lem:quadric_relations}
  In this situation, the $(n)$-ruling $D$ gives the following $k-(n+1)$
  quadratic relations in $\Cox(S_r)$: \[q_i := \sum_{j=1}^k a_{j,i} \cdot
  \xi(D_j)\cdot\xi(D'_j) = 0 \qquad\text{for $i \in \{1, \dots, k-(n+1)\}.$}\]
\end{lemma}

We will describe the $(n)$-rulings in more detail in the subsequent sections.

Let $J_r$ be the ideal in $R_r$ which is generated by the $(n)$-rulings (where
$n=1$ for $r \le 6$, $n \in \{1,2\}$ for $r=7$, and $n\in \{1,2,3\}$ for
$r=8$).

\textbf{The proof of Theorem~\ref{thm:relations}.} For $r \in \{4,5,6\}$, this
is \cite[Theorem~4.9]{MR2029863}. For $r \in \{7,8\}$, we use a refinement of
its proof.

Let $Z_r = \Spec(R_r/\rad(J_r)) \subset \Spec(R_r)$. We want to prove
that $Z_r$ equals $\AA(S_r) \subset \Spec(R_r)$, where
$\AA(S_r):=\Spec(\Cox(S_r))$. Obviously, $0 \in \Spec(R_r)$ is
contained in both $Z_r$ and $\AA(S_r)$. Its complement $\Spec(R_r)
\setminus \{0\}$ is covered by the open sets \[U_D := \{\xi(D) \ne
0\},\qquad\text{where $D \in \DD_r$.}\] In the case $r = 8$, we will
show that it suffices to consider the sets $U_D$ for $D \in \DD_8
\setminus \{K_1,K_2\}$.

We want to show \[Z_r \cap U_D \cong Z_{r-1} \times (\AA^1\setminus
\{0\}).\] Note that we can identify the negative curves $\DD_{r-1}$ of
$S_{r-1}$ with the subset $\DD'_r$ of $\DD_r$ containing the negative
curves which do not intersect $D$. We define
\[\begin{array}{cccc}
  \psi: & Z_r \cap U_D & \to & Z_{r-1} \times (\AA^1 \setminus \{0\})\\
  & (\xi(D') \mid D' \in \DD_r) & \mapsto & 
  ((\xi(D') \mid D' \in \DD_{r-1}), \xi_D)
\end{array}.\]

For $r \in \{7,8\}$, we will prove:

\begin{lemma}\label{lem:dependence}
  Every $\xi(D'')$ for $D'' \in \DD_r$ intersecting $D$ is determined by
  \[\xi(D)\qquad\text{and}\qquad\{\xi(D') \mid D' \in \DD_r
  \quad\text{with}\quad (D',D) = 0\},\]
  provided that $\xi(D) \ne 0$ and using the relations generating $J_r$.
\end{lemma}

By the proof of \cite[Prop. 4.4]{MR2029863}, 
\[\AA(S_r) \cap U_D \cong \AA(S_{r-1}) \times (\AA^1 \setminus \{0\}).\]
By induction, $Z_{r-1} = \AA(S_{r-1})$. Therefore, $Z_r \cap U_D =
\AA(S_r) \cap U_D$ for every negative curve $D$, which implies $Z_r =
\AA(S_r)$, completing the proof of Theorem~\ref{thm:relations} once
Lemma~\ref{lem:dependence} is proved.

\

\textbf{Proof of Theorem~\ref{thm:radical}.} We want to show that
the ideal $J_r$ is radical.

\begin{lemma}\label{lem:hilbert}
  The Hilbert polynomial of $S_r/J_r$ has degree $r+2$.
\end{lemma}

For $r=5$, this was proved by Batyrev. We will prove it for $r \in
\{6,7\}$.

\begin{remark}
  The problem of calculating the Hilbert polynomial of $J_8$ seems 
  out of reach of the current computer algebra packages. It is the
  only step missing in the proof of Theorem~\ref{thm:radical} for Del
  Pezzo surfaces of degree~$1$.
\end{remark}

Under the condition of the proof of Lemma~\ref{lem:hilbert}, the depth
of $R_r/J_r$ is $r+3$. As $\Spec(R_r/J_r)$ is irreducible by
\cite{MR2029863}, and $\Cox(S_r) = \Spec(R_r/\rad(J_r))$ by
Theorem~\ref{thm:relations}, the $R_r$-module $R_r/J_r$ is
Cohen-Macaulay. Therefore, we need to check the following claim in
order to prove that the ideal $J_r$ is radical:

\begin{lemma}\label{lem:smooth_point}
  $R_r/J_r$ has a smooth point.
\end{lemma}

\section{Degree 3}

We consider the case $r=6$, i.e., smooth cubic surfaces. The set
$\DD_6$ of negative curves on $S_6$ consists of the following 27
divisors:
\begin{itemize}
\item exceptional divisors $E_1, \dots, E_6$, preimages of
  $p_1, \dots, p_6 \in \Ptwo$,
\item transforms $m_{i,j} = H-E_i-E_j$ of the 15 lines $m'_{i,j}$ through
  the points $p_i,p_j$ ($i\ne j \in \{1, \dots, 6\}$), and
\item transforms $Q_k = H - (E_1+\dots+E_6)+E_k$ of the six conics $Q'_k$
  through all of the blown-up points except $p_k$.
\end{itemize}
With respect to the anticanonical embedding $S_6 \inj \Pthree$, the
negative curves are the 27 lines.

Together with information from Section~\ref{sec:generalities}, it
is straightforward to derive: 

\begin{lemma}\label{lem:dynkin_degree3}
  The extended Dynkin diagram of negative curves has the following structure:
  \begin{enumerate}
  \item It has 27 vertices corresponding to the 27 lines $E_i, m_{i,j},
    Q_i$. Each of them has self-intersection number $-1$.
  \item Every line intersects exactly 10 other lines: $E_i$ intersects
    $m_{i,j}$ and $Q_j$ (for $j \ne i$); $m_{i,j}$ intersects
    $E_i,E_j,Q_i,Q_j$ and $m_{k,l}$ (for $\{i,j\}\cap\{k,l\} = \emptyset$);
    $Q_i$ intersects $m_{i,j}$ and $E_j$ (for $j \ne i$).  Correspondingly,
    there are 135 edges in the Dynkin diagram.
  \item\label{item:dynkin_degree3_triangles} There are 45
    \emph{triangles}, i.e., triples of lines which intersect pairwise:
    30 triples $E_i, m_{i,j}, Q_j$ and 15 triples $m_{i_1,j_1},
    m_{i_2,j_2}, m_{i_3,j_3}$ where $\{i_1,j_1,i_2,j_2,i_3,j_3\} =
    \{1, \dots, 6\}$. This corresponds to 45 triangles in the Dynkin
    diagram, where each edge is contained in exactly one of the
    triangles, and each vertex belongs to exactly five triangles.
  \end{enumerate}
\end{lemma}

\begin{lemma}
  The 27 rulings of $S_6$ are given by $\anti 6 - D$ for $D \in
  \DD_6$. Two negative curves $D', D''$ fulfill $D'+D'' = \anti 6 - D$
  if and only if $D, D', D''$ form a triangle in the sense of
  Lemma~\ref{lem:dynkin_degree3}(\ref{item:dynkin_degree3_triangles}).
  There are five such pairs for any given $D$.
\end{lemma}

\begin{proof}
  We can check directly that $D+D'+D'' = \anti 6$ if $D,D',D''$ form a
  triangle. Therefore, $\anti 6 - D$ is a ruling, and as any $D$
  is contained in exactly five triangles, it can be expressed in five
  corresponding ways as $D'+D''$.
  
  On the other hand, by Table~\ref{tab:relations}, the total number of
  rulings is 27, and each ruling can be expressed in exactly five ways
  as the sum of two negative curves.
\end{proof}

Let $D$ be one of the 27 lines of $S_6$, and consider the projection
$\psi_D: S_6 \rto \Ptwo$ from $D$. Then \[\psi_D^*(\OO_\Ptwo(1)) =
-K_{S_6}-D=
\begin{cases}
  H-E_i, &D = Q_i,\\
  2H-(E_1+\dots+E_6)+E_i+E_j, &D = m_{i,j},\\
  3H-(E_1+\dots+E_6)-E_i, &D = E_i.
\end{cases}
\]
These are exactly the rulings.

A generating set of $\Cox(S_6)$ is given by section $\eta_i,
\mu_{i,j}, \lambda_i$ corresponding to the 27 lines $E_i, m_{i,j},
Q_i$, respectively. Let \[R_6 := \KK[\eta_i, \mu_{i,j}, \lambda_i].\]

The quadratic monomials in $\Gamma(S_6,\anti 6 - D)$ corresponding to the
five ways to express $\anti 6-D$ as the sum of the negative curves are
\begin{itemize}
\item $\mu_{i,j}\eta_j$ if $D = Q_i$
\item $\eta_i\lambda_j, \eta_j\lambda_i, \mu_{k_1,k_2}\mu_{k_3,k_4}$
  if $D = \mu_{i,j}$ (with $\{i,j,k_1, \dots, k_4\} = \{1, \dots,
  6\}$)
\item $\mu_{i,j}\lambda_j$ if $D = E_i$
\end{itemize}

In order to calculate the 81 relations in $J_6$ explicitly as
described in Lemma~\ref{lem:quadric_relations}, we use the coordinates
of \eqref{eq:four_points} for $p_1, \dots, p_4$, and
\[p_5 = (1:a:b), \qquad p_6 = (1:c:d).\]
We write \[E:=(b-1)(c-1)-(a-1)(d-1)\qquad\text{and}\qquad F:=bc-ad\]
for simplicity. The three relations corresponding to a line $D$
will be denoted by $\qa D,\qb D,\qc D$.
\begin{align*}
\qa{Q_1}&=-\eta_2\mu_{1,2}-\eta_3\mu_{1,3}+\eta_4\mu_{1,4}\\
\qb{Q_1}&=-a\eta_2\mu_{1,2}-b\eta_3\mu_{1,3}+\eta_5\mu_{1,5}\\
\qb{Q_1}&=-c\eta_2\mu_{1,2}-d\eta_3\mu_{1,3}+\eta_6\mu_{1,6}
\displaybreak[0]\\[\baselineskip]
\qa{Q_2}&=\eta_1\mu_{1,2}-\eta_3\mu_{2,3}+\eta_4\mu_{2,4}\\
\qb{Q_2}&=\eta_1\mu_{1,2}-b\eta_3\mu_{2,3}+\eta_5\mu_{2,5}\\
\qc{Q_2}&=\eta_1\mu_{1,2}-d\eta_3\mu_{2,3}+\eta_6\mu_{2,6}
\displaybreak[0]\\[\baselineskip]
\qa{Q_3}&=\eta_1\mu_{1,3}+\eta_2\mu_{2,3}+\eta_4\mu_{3,4}\\
\qb{Q_3}&=\eta_1\mu_{1,3}+a\eta_2\mu_{2,3}+\eta_5\mu_{3,5}\\
\qc{Q_3}&=\eta_1\mu_{1,3}+c\eta_2\mu_{2,3}+\eta_6\mu_{3,6}
\displaybreak[0]\\[\baselineskip]
\qa{Q_4}&=\eta_1\mu_{1,4}+\eta_2\mu_{2,4}+\eta_3\mu_{3,4}\\
\qb{Q_4}&=(1-b)\eta_1\mu_{1,4}+(a-b)\eta_2\mu_{2,4}+\eta_5\mu_{4,5}\\
\qc{Q_4}&=(1-d)\eta_1\mu_{1,4}+(c-d)\eta_2\mu_{2,4}+\eta_6\mu_{4,6}
\displaybreak[0]\\[\baselineskip]
\qa{Q_5}&=1/b\eta_1\mu_{1,5}+a/b\eta_2\mu_{2,5}+\eta_3\mu_{3,5}\\
\qb{Q_5}&=(1-b)/b\eta_1\mu_{1,5}+(a-b)/b\eta_2\mu_{2,5}+\eta_4\mu_{4,5}\\
\qc{Q_5}&=(b-d)/b\eta_1\mu_{1,5}+F/b\eta_2\mu_{2,5}+\eta_6\mu_{5,6}
\displaybreak[0]\\[\baselineskip]
\qa{Q_6}&=1/d\eta_1\mu_{1,6}+c/d\eta_2\mu_{2,6}+\eta_3\mu_{3,6}\\
\qb{Q_6}&=(1-d)/d\eta_1\mu_{1,6}+(c-d)/d\eta_2\mu_{2,6}+\eta_4\mu_{4,6}\\
\qc{Q_6}&=(b-d)/d\eta_1\mu_{1,6}+F/d\eta_2\mu_{2,6}+\eta_5\mu_{5,6}
\displaybreak[0]\\[\baselineskip]
\qa{m_{1,2}}&=\mu_{4,5}\mu_{3,6}-\mu_{3,5}\mu_{4,6}+\mu_{3,4}\mu_{5,6}\\
\qb{m_{1,2}}&=(b-d)\mu_{3,5}\mu_{4,6}+(d-1)\mu_{3,4}\mu_{5,6}+\eta_2\lambda_1\\
\qc{m_{1,2}}&=F\mu_{3,5}\mu_{4,6}+a(d-c)\mu_{3,4}\mu_{5,6}+\eta_1\lambda_2
\displaybreak[0]\\[\baselineskip]
\qa{m_{1,3}}&=\mu_{4,5}\mu_{2,6}-\mu_{2,5}\mu_{4,6}+\mu_{2,4}\mu_{5,6}\\
\qb{m_{1,3}}&=(c-a)\mu_{2,5}\mu_{4,6}+(1-c)\mu_{2,4}\mu_{5,6}+\eta_3\lambda_1\\
\qc{m_{1,3}}&=-F\mu_{2,5}\mu_{4,6}+b(c-d)\mu_{2,4}\mu_{5,6}+\eta_1\lambda_3
\displaybreak[0]\\[\baselineskip]
\qa{m_{2,3}}&=\mu_{4,5}\mu_{1,6}-\mu_{1,5}\mu_{4,6}+\mu_{1,4}\mu_{5,6}\\
\qb{m_{2,3}}&=(a-c)\mu_{1,5}\mu_{4,6}+a(c-1)\mu_{1,4}\mu_{5,6}+\eta_3\lambda_2\\
\qc{m_{2,3}}&=(b-d)\mu_{1,5}\mu_{4,6}+b(d-1)\mu_{1,4}\mu_{5,6}+\eta_2\lambda_3
\displaybreak[0]\\[\baselineskip]
\qa{m_{1,4}}&=\mu_{3,5}\mu_{2,6}-\mu_{2,5}\mu_{3,6}+\mu_{2,3}\mu_{5,6}\\
\qb{m_{1,4}}&=-E\mu_{2,5}\mu_{3,6}+(b-1)(c-1)\mu_{2,3}\mu_{5,6}+\eta_4\lambda_1\\
\qc{m_{1,4}}&=-F\mu_{2,5}\mu_{3,6}+bc\mu_{2,3}\mu_{5,6}+\eta_1\lambda_4
\displaybreak[0]\\[\baselineskip]
\qa{m_{2,4}}&=\mu_{3,5}\mu_{1,6}-\mu_{1,5}\mu_{3,6}+\mu_{1,3}\mu_{5,6}\\
\qb{m_{2,4}}&=E\mu_{1,5}\mu_{3,6}+(a-b)(c-1)\mu_{1,3}\mu_{5,6}+\eta_4\lambda_2\\
\qc{m_{2,4}}&=(b-d)\mu_{1,5}\mu_{3,6}-b\mu_{1,3}\mu_{5,6}+\eta_2\lambda_4
\displaybreak[0]\\[\baselineskip]
\qa{m_{3,4}}&=\mu_{2,5}\mu_{1,6}-\mu_{1,5}\mu_{2,6}+\mu_{1,2}\mu_{5,6}\\
\qb{m_{3,4}}&=-E\mu_{1,5}\mu_{2,6}+(a-b)(1-d)\mu_{1,2}\mu_{5,6}+\eta_4\lambda_3\\
\qc{m_{3,4}}&=(c-a)\mu_{1,5}\mu_{2,6}+a\mu_{1,2}\mu_{5,6}+\eta_3\lambda_4
\displaybreak[0]\\[\baselineskip]
\qa{m_{1,5}}&=\mu_{3,4}\mu_{2,6}-\mu_{2,4}\mu_{3,6}+\mu_{2,3}\mu_{4,6}\\
\qb{m_{1,5}}&=-E\mu_{2,4}\mu_{3,6}+(a-c)(1-b)\mu_{2,3}\mu_{4,6}+\eta_5\lambda_1\\
\qc{m_{1,5}}&=(d-c)\mu_{2,4}\mu_{3,6}+c\mu_{2,3}\mu_{4,6}+\eta_1\lambda_5
\displaybreak[0]\\[\baselineskip]
\qa{m_{2,5}}&=\mu_{3,4}\mu_{1,6}-\mu_{1,4}\mu_{3,6}+\mu_{1,3}\mu_{4,6}\\
\qb{m_{2,5}}&=aE\mu_{1,4}\mu_{3,6}+(a-b)(c-a)\mu_{1,3}\mu_{4,6}+\eta_5\lambda_2\\
\qc{m_{2,5}}&=(1-d)\mu_{1,4}\mu_{3,6}-\mu_{1,3}\mu_{4,6}+\eta_2\lambda_5
\displaybreak[0]\\[\baselineskip]
\qa{m_{3,5}}&=\mu_{2,4}\mu_{1,6}-\mu_{1,4}\mu_{2,6}+\mu_{1,2}\mu_{4,6}\\
\qb{m_{3,5}}&=-bE\mu_{1,4}\mu_{2,6}+(a-b)(b-d)\mu_{1,2}\mu_{4,6}+\eta_5\lambda_3\\
\qc{m_{3,5}}&=(c-1)\mu_{1,4}\mu_{2,6}+\mu_{1,2}\mu_{4,6}+\eta_3\lambda_5
\displaybreak[0]\\[\baselineskip]
\qa{m_{4,5}}&=\mu_{2,3}\mu_{1,6}-\mu_{1,3}\mu_{2,6}+\mu_{1,2}\mu_{3,6}\\
\qb{m_{4,5}}&=b(c-a)\mu_{1,3}\mu_{2,6}+a(b-d)\mu_{1,2}\mu_{3,6}+\eta_5\lambda_4\\
\qc{m_{4,5}}&=(c-1)\mu_{1,3}\mu_{2,6}+(1-d)\mu_{1,2}\mu_{3,6}+\eta_4\lambda_5
\displaybreak[0]\\[\baselineskip]
\qa{m_{1,6}}&=\mu_{3,4}\mu_{2,5}-\mu_{2,4}\mu_{3,5}+\mu_{2,3}\mu_{4,5}\\
\qb{m_{1,6}}&=-E\mu_{2,4}\mu_{3,5}+(a-c)(1-d)\mu_{2,3}\mu_{4,5}+\eta_6\lambda_1\\
\qc{m_{1,6}}&=(b-a)\mu_{2,4}\mu_{3,5}+a\mu_{2,3}\mu_{4,5}+\eta_1\lambda_6
\displaybreak[0]\\[\baselineskip]
\qa{m_{2,6}}&=\mu_{3,4}\mu_{1,5}-\mu_{1,4}\mu_{3,5}+\mu_{1,3}\mu_{4,5}\\
\qb{m_{2,6}}&=cE\mu_{1,4}\mu_{3,5}+(a-c)(d-c)\mu_{1,3}\mu_{4,5}+\eta_6\lambda_2\\
\qc{m_{2,6}}&=(1-b)\mu_{1,4}\mu_{3,5}-\mu_{1,3}\mu_{4,5}+\eta_2\lambda_6
\displaybreak[0]\\[\baselineskip]
\qa{m_{3,6}}&=\mu_{2,4}\mu_{1,5}-\mu_{1,4}\mu_{2,5}+\mu_{1,2}\mu_{4,5}\\
\qb{m_{3,6}}&=-dE\mu_{1,4}\mu_{2,5}+(d-b)(d-c)\mu_{1,2}\mu_{4,5}+\eta_6\lambda_3\\
\qc{m_{3,6}}&=(a-1)\mu_{1,4}\mu_{2,5}+\mu_{1,2}\mu_{4,5}+\eta_3\lambda_6
\displaybreak[0]\\[\baselineskip]
\qa{m_{4,6}}&=\mu_{2,3}\mu_{1,5}-\mu_{1,3}\mu_{2,5}+\mu_{1,2}\mu_{3,5}\\
\qb{m_{4,6}}&=d(c-a)\mu_{1,3}\mu_{2,5}+c(b-d)\mu_{1,2}\mu_{3,5}+\eta_6\lambda_4\\
\qc{m_{4,6}}&=(a-1)\mu_{1,3}\mu_{2,5}+(1-b)\mu_{1,2}\mu_{3,5}+\eta_4\lambda_6
\displaybreak[0]\\[\baselineskip]
\qa{m_{5,6}}&=\mu_{2,3}\mu_{1,4}-\mu_{1,3}\mu_{2,4}+\mu_{1,2}\mu_{3,4}\\
\qb{m_{5,6}}&=d(c-1)\mu_{1,3}\mu_{2,4}+c(1-d)\mu_{1,2}\mu_{3,4}+\eta_6\lambda_5\\
\qc{m_{5,6}}&=b(a-1)\mu_{1,3}\mu_{2,4}+a(1-b)\mu_{1,2}\mu_{3,4}+\eta_5\lambda_6
\displaybreak[0]\\[\baselineskip]
\qa{E_1}&=(d-b)/E\mu_{1,2}\lambda_2+(c-a)/E\mu_{1,3}\lambda_3+\mu_{1,4}\lambda_4\\
\qb{E_1}&=(d-1)/E\mu_{1,2}\lambda_2+(c-1)/E\mu_{1,3}\lambda_3+\mu_{1,5}\lambda_5\\
\qc{E_1}&=(b-1)/E\mu_{1,2}\lambda_2+(a-1)/E\mu_{1,3}\lambda_3+\mu_{1,6}\lambda_6
\displaybreak[0]\\[\baselineskip]
\qa{E_2}&=F/E\mu_{1,2}\lambda_1+(c-a)/E\mu_{2,3}\lambda_3+\mu_{2,4}\lambda_4\\
\qb{E_2}&=(c-d)/E\mu_{1,2}\lambda_1+(c-1)/E\mu_{2,3}\lambda_3+\mu_{2,5}\lambda_5\\
\qc{E_2}&=(a-b)/E\mu_{1,2}\lambda_1+(a-1)/E\mu_{2,3}\lambda_3+\mu_{2,6}\lambda_6
\displaybreak[0]\\[\baselineskip]
\qa{E_3}&=F/E\mu_{1,3}\lambda_1+(b-d)/E\mu_{2,3}\lambda_2+\mu_{3,4}\lambda_4\\
\qb{E_3}&=(c-d)/E\mu_{1,3}\lambda_1+(1-d)/E\mu_{2,3}\lambda_2+\mu_{3,5}\lambda_5\\
\qc{E_3}&=(a-b)/E\mu_{1,3}\lambda_1+(1-b)/E\mu_{2,3}\lambda_2+\mu_{3,6}\lambda_6
\displaybreak[0]\\[\baselineskip]
\qa{E_4}&=F/(a-c)\mu_{1,4}\lambda_1+(b-d)/(a-c)\mu_{2,4}\lambda_2+\mu_{3,4}\lambda_3\\
\qb{E_4}&=c/(a-c)\mu_{1,4}\lambda_1+1/(a-c)\mu_{2,4}\lambda_2+\mu_{4,5}\lambda_5\\
\qc{E_4}&=a/(a-c)\mu_{1,4}\lambda_1+1/(a-c)\mu_{2,4}\lambda_2+\mu_{4,6}\lambda_6
\displaybreak[0]\\[\baselineskip]
\qa{E_5}&=(d-c)/(c-1)\mu_{1,5}\lambda_1+(d-1)/(c-1)\mu_{2,5}\lambda_2+\mu_{3,5}\lambda_3\\
\qb{E_5}&=-c/(c-1)\mu_{1,5}\lambda_1-1/(c-1)\mu_{2,5}\lambda_2+\mu_{4,5}\lambda_4\\
\qc{E_5}&=-1/(c-1)\mu_{1,5}\lambda_1-1/(c-1)\mu_{2,5}\lambda_2+\mu_{5,6}\lambda_6
\displaybreak[0]\\[\baselineskip]
\qa{E_6}&=(b-a)/(a-1)\mu_{1,6}\lambda_1+(b-1)/(a-1)\mu_{2,6}\lambda_2+\mu_{3,6}\lambda_3\\
\qb{E_6}&=-a/(a-1)\mu_{1,6}\lambda_1-1/(a-1)\mu_{2,6}\lambda_2+\mu_{4,6}\lambda_4\\
\qc{E_6}&=-1/(a-1)\mu_{1,6}\lambda_1-1/(a-1)\mu_{2,6}\lambda_2+\mu_{5,6}\lambda_5
\end{align*}

\bigskip

\textbf{Proof of Lemma~\ref{lem:hilbert}.} We calculate the
Hilbert polynomial of $R_6/J_6$ over the field of fractions of the
polynomial ring $\QQ[a,b,c,d]$ using \texttt{Magma}:
\begin{multline*}
  h(t) = \frac {1}{8!} (372t^8 + 4464t^7 + 25200t^6 + 86184t^5 + 193788t^4\\
  + 291816t^3 + 284640t^2 + 161856t + 40320).
\end{multline*}
Its degree is $r+2=8$ as required.

\

\textbf{Proof of Lemma~\ref{lem:smooth_point}.} The point $p$ with
coordinates
\[(\eta_5,\eta_6,\mu_{1,2},\mu_{1,4},\mu_{2,3},\mu_{3,4},\lambda_5,\lambda_6)
= (c(d-1),a(b-1),1,-1,1,1,1,1),\] and all other coordinates zero is a
smooth point of $R_6/J_6$. Indeed, we check that $p$ fulfills all the
relations generating $J_r$ (which is obvious for all of them except
$\qa{m_{5,6}}, \qb{m_{5,6}}, \qc{m_{5,6}}$), and we calculate directly
that the $81 \times 27$ Jacobian matrix has full rank 18 at this
point.

\section{Degree 2}

Let $S_7$ be a smooth Del Pezzo surface of degree $d=2$, i.e., the
blow-up of $\Ptwo$ in $r=7$ points.
The set $\DD_7$ contains 56 negatives curves which are the transforms of
the following curves in $\Ptwo$:
\begin{itemize}
\item blow-ups $E_1, \dots, E_7$ of $p_1,\dots, p_7$;
\item 21 lines $m_{i,j}'$ through $p_i,p_j$, where
  \[m_{i,j} = H-E_i-E_j;\]
\item 21 conics $Q_{i,j}'$ through five of the seven points,
  missing $p_i,p_j$, where
  \[Q_{i,j} = 2H-(E_1+\dots+E_7)+E_i+E_j;\]
\item 7 singular cubics $C_i'$ through all seven points, where
  $p_i$ is a double point, and
  \[C_i = 3H-(E_1+\dots+E_7)-E_i.\]
\end{itemize}

The Cox ring $\Cox(S_7)$ is generated by the sections $\eta_i,
\mu_{i,j}, \nu_{i,j}, \lambda_i$ corresponding to the 56 negative
curves $E_i, m_{i,j}, Q_{i,j}, C_i$, respectively. Let \[R_7 :=
\KK[\eta_i, \mu_{i,j}, \nu_{i,j}, \lambda_i]\] be the polynomial ring
in 56 generators.

Consider the ideal $I_7 \subset R_7$ generated by the quadratic
relations corresponding to rulings. In view of
Lemma~\ref{lem:quadric_relations}, we need to know the six different
ways to write each of the 126 rulings as a sum of two negative curves
in order to describe $I_7$ explicitly. Here, we do not write the
resulting 504 relations down because of the length of this list.

\begin{lemma}
  Each of the 126 rulings can be written in six ways as a sum of two negative
  curves:
  \begin{enumerate}
  \item For the seven rulings $H-E_i$:
    \[\{E_j + m_{i,j} \mid j \ne i\}.\]
  \item For the 35 rulings $2H-(E_1+\dots+E_7)+E_i+E_j+E_k$:
   {\small \[\{E_i+Q_{j,k}, E_j+Q_{i,k}, E_k+Q_{i,j}, m_{l_1,l_2}+m_{l_3,l_4} \mid
    \{ i,j,k,l_1,l_2,l_3,l_4\} = \{1, \dots, 7\}\}\]}
  \item For the 42 rulings $3H-(E_1+\dots+E_7)+E_i-E_j$:
    \[\{E_i+C_j, Q_{i,k}+m_{j,k} \mid k \ne i,j\}.\]
  \item For the 35 rulings $4H-(E_1+\dots+E_7)-E_i-E_j-E_k$:
   {\small \[\{C_i+m_{j,k}, C_j+m_{i,k}, C_k+m_{i,j}, Q_{l_1,l_2}+Q_{l_3,l_4} \mid
    \{i,j,k,l_1,l_2,l_3,l_4\} = \{1, \dots, 7\}\}\]}
  \item For the seven rulings $5H-2(E_1+\dots+E_7)+E_i$:
    \[\{C_j + Q_{i,j} \mid j \ne i\}.\]
  \end{enumerate}
\end{lemma}

However, we have more quadratic relations in $\Cox(S_7)$: Note that
the point $q$, with $\eta_1=\lambda_1=1$ and other coordinates
zero, satisfies the 504 relations. Indeed, $(E_1,C_1) = 2$, but all
quadratic monomials which occur in the relations correspond to pairs
of divisors whose intersection number is $1$. Therefore, all these
monomials and all the relations vanish in $q$. On the other hand, we
check that the $504 \times 56$ Jacobian matrix has rank 54 in this
point, which means that $q$ is contained in a component of the variety
defined by $I_7$ which has dimension~$2$. As $\AA(S_7)$ is irreducible
of dimension $10$, we must find other relations to exclude such
components.

As $E_1+C_1 = \anti 7$, we look for more relations in degree $\anti 7$ of
$\Cox(S_7)$: We check that in this degree, we have exactly 28 monomials:
\[\{\eta_i\lambda_i \mid 1 \le i \le 7\}\cup\{\mu_{j,k}\nu_{j,k} \mid 
1 \le j < k \le 7\},\] corresponding to $\anti 7 = E_i+C_i =
m_{j,k}+Q_{j,k}$. As $\dim \Gamma(S_7, \anti 7) = 3$, and as none of the
relations coming from rulings induces a relation in this degree, we
obtain 25 independent relations. Note that $\anti 7$ is the unique
$(2)$-ruling of $S_7$.

We can calculate the relations explicitly as they correspond to the
relations between the polynomials $f_{E_i}\cdot f_{C_i}$ and
$f_{m_{i,j}} \cdot f_{Q_{i,j}}$, which are homogeneous of degree 3, as
described in Lemma~\ref{lem:quadric_relations}.

Let $J_7$ be the ideal generated by these $529$ relations.

\

\textbf{Proof of Lemma~\ref{lem:dependence}.} In order to show
that $\Cox(S_7)$ is described by $\rad(J_7)$, we must prove
Lemma~\ref{lem:dependence} in the case $r=7$.

For any $D \in \DD_7$, consider a coordinate $\xi(D')$ where $(D,D')
= 1$. This is determined by the ruling $D+D'$. Indeed, this ruling
induces a relation of the form \[\xi(D)\xi(D') = \sum
a_i\xi(D_i)\xi(D_i'),\] where $D_i+D_i' = D+D'$. Therefore,
\[(D,D_i+D_i') = (D,D+D') = (D,D)+(D,D') =-1+1= 0,\] which implies 
$(D,D_i)=(D,D_i') = 0$ since the only negative curve intersecting $D$
negatively is $D$ itself. Since $\xi(D) \ne 0$, the only unknown
variable $\xi(D')$ is determined by this relation.

Furthermore, there is exactly one coordinate $\xi(D'')$ where $(D,D'')
= 2$. The unique $(2)$-ruling $D+D'' = \anti 7$ induces a relation of
the form
\[\xi(D)\xi(D'') = \sum a_i\xi(D_i)\xi(D'_i),\] where $\xi(D'')$ is 
the only unknown variable.

\

\textbf{Proof of Lemma~\ref{lem:hilbert}.} In a special case, we
can calculate the Hilbert polynomial:

\begin{example}
  Over the field $\FF_{101}$ with $p_1, \dots, p_4$ as
  in~\eqref{eq:four_points} and \[p_5 = (1:2:3), \quad p_6 = (1:5:7), \quad
  p_7 = (1:13:17)\] in general position, we can use \texttt{Macaulay2} to
  calculate the Hilbert polynomial of $J_7$ as
  \begin{multline}\label{hilbert_deg2}
    h(t) = \frac 1 {9!}\cdot (9504 t^9 + 85536 t^8 + 412992 t^7 + 1294272 t^6 
    + 2860704 t^5\\ + 4554144 t^4 + 5125248 t^3 + 3863808 t^2
    + 1752192 t + 362880).
  \end{multline}
\end{example}

The Hilbert polynomial does not depend on the choice of the field or
the points. Therefore, $h(t)$ is the Hilbert polynomial of $R_7/J_7$.
Its degree is $r+2=9$.

\section{Degree 1}\label{sec:degree1}

In this section, we consider blow-ups of $\Ptwo$ in $r=8$ points in
general position, i.e., Del Pezzo surfaces $S_8$ of degree $1$.

The set $\DD_8$ contains the transforms of the following 242 curves:
\begin{itemize}
\item Blow-ups $E_1, \dots, E_8$ of $p_1, \dots, p_8$;
\item $28$ lines $m_{i,j}'$ through $p_i, \, p_j$:
  \[m_{i,j} = H-E_i-E_j;\]
\item $56$ conics $Q_{i,j,k}'$ through $5$ points, missing
  $p_i,\,p_j,\,p_k$: \[Q_{i,j,k} = 2H-(E_1+\dots+E_8)+E_i+E_j+E_k;\]
\item $56$ cubics $C_{i,j}'$ through $7$ points missing $p_j$,
  where $p_i$ is a double point:
  \[C_{i,j} = 3H-(E_1+\dots+E_8)-E_i+E_j;\]
\item $56$ quartics $V_{i,j,k}'$ through all points, where
  $p_i,\,p_j,\,p_k$ are double points:
  \[V_{i,j,k} = 4H-(E_1+\dots+E_8)-(E_i+E_j+E_k);\]
\item $28$ quintics $F_{i,j}'$ through all points, where $p_i,\,p_j$ are
  simple points and the other six are double points:
  \[F_{i,j} = 5H-2(E_1+\dots+E_8)+E_i+E_j;\]
\item $8$ sextics $T_i'$, where $p_i$ a triple point and the other
  seven points are double points:
  \[T_i = 6H-2(E_1+\dots+E_8)-E_i;\]
\item two independent cubics $K_1', K_2'$ through the eight points:
  \[[K_1]=[K_2]=\anti 8 = 3-(E_1+\dots+E_8).\]
\end{itemize}
The Cox ring of $S_8$ is generated by the 242 sections
\[\eta_i,\, \mu_{i,j},\, \nu_{i,j,k},\,\lambda_{i,j},\,
\phi_{i,j,k},\,\psi_{i,j},\,\sigma_i,\,\kappa_i\] of $E_i,\,
m_{i,j},\,Q_{i,j,k},\,C_{i,j},\,V_{i,j,k},\,F_{i,j},\,T_i, K_i$,
respectively.

\begin{lemma}
  Each of the 2160 rulings can be expressed in the following seven
  ways as a sum of two negative curves:
  \begin{itemize}
  \item $8$ rulings of the form $H-E_i$: \[\{E_j+m_{i,j} \mid j \ne
    i\}.\]
  \item $\binom 8 4 = 70$ rulings of the form $2H-(E_i+E_j+E_k+E_l)$:
    \[\Bigg\{
    \begin{aligned}
      &m_{i,j}+m_{k,l},m_{i,k}+m_{j,l},\\
      &m_{i,l}+m_{j,k},E_a+Q_{b,c,d}
    \end{aligned}
    \Bigg| \{a,b,c,d,i,j,k,l\} = \{1,\dots,8\}\Bigg\}.\]
  \item $8 \cdot \binom 7 2 = 168$ rulings of the form
    $3H-(E_1+\dots+E_8)-E_i+E_j+E_k$:
    \[\{E_j+C_{i,k},E_k+C_{i,j},m_{i,l}+Q_{j,k,l} \mid l \notin \{i,j,k\}\}.\]
  \item $8 \cdot \binom 7 3 = 280$ rulings 
    $4H-(E_1+\dots+E_8)+E_i-(E_j+E_k+E_l)$:
    \[\Bigg\{
    \begin{aligned}
      &E_i+V_{j,k,l},Q_{i,a,b}+Q_{i,c,d},\\
      &C_{j,i}+m_{k,l},C_{k,i}+m_{j,l},C_{l,i}+m_{j,k}
    \end{aligned}
    \Bigg| \{a,b,c,d,i,j,k,l\} = \{1,\dots,8\}\Bigg\},\]
    and $8$ rulings of the form $4H-(E_1+\dots+E_8)-2E_i$:
    \[\{m_{i,j}+C_{i,j} \mid j \ne i\}.\]
  \item $8\cdot 7 = 56$ rulings of the form
    $5H-2(E_1+\dots+E_8)+2E_i+E_j$:
    \[\{E_i+F_{i,j},C_{k,i}+Q_{i,j,k} \mid k \notin \{i,j\}\},\]
    and $8\cdot \binom 7 3 = 280$ rulings 
    $5H-(E_1+\dots+E_8)-2E_i-(E_j+E_k+E_l)$:
    \[\Bigg\{
    \begin{aligned}
      &m_{i,j}+V_{i,k,l},m_{i,k}+V_{i,j,l},\\
      &m_{i,l}+V_{i,j,k}, C_{i,a}+Q_{b,c,d}
    \end{aligned}
    \Bigg| \{a,b,c,d,i,j,k,l\} = \{1,\dots,8\}\Bigg\}.\]
  \item $\binom 8 2 \cdot \binom 6 2 = 420$ rulings
    $6H-2(E_1+\dots+E_8)-(E_i+E_j)+E_k+E_l$:
    \[\{m_{i,j}+F_{k,l},V_{i,j,m}+Q_{k,l,m},C_{i,k}+C_{j,l},C_{i,l}+C_{j,k} 
    \mid m \notin \{i,j,k,l\}\}.\]
  \item $8 \cdot 7 = 56$ rulings of the form $7H-2(E_1+\dots+E_8)-2E_i-E_j$:
    \[\{m_{i,j}+T_i,C_{i,k}+V_{i,j,k} \mid k \notin \{i,j\}\},\]
    and $8\cdot \binom 7 3 = 280$ rulings 
    $7H-3(E_1+\dots+E_8)+2E_i+E_j+E_k+E_l$:
    \[\Bigg\{
    \begin{aligned}
      &F_{i,j}+Q_{i,k,l},F_{i,k}+Q_{i,j,l},\\
      &F_{i,l}+Q_{i,j,k},C_{a,i}+V_{b,c,d}
  \end{aligned}
 \Bigg| \{a,b,c,d,i,j,k,l\} = \{1,\dots,8\}\Bigg\}.\]
  \item $8\cdot \binom 7 3 = 280$ rulings
    $8H-3(E_1+\dots+E_8)-E_i+E_j+E_k+E_l$:
    \[\Bigg\{
    \begin{aligned}
      &C_{i,j}+F_{k,l},C_{i,k}+F_{j,l},C_{i,l}+F_{j,k},\\
      &T_i+Q_{j,k,l},V_{i,a,b}+V_{i,c,d}
  \end{aligned}
  \Bigg|\{a,b,c,d,i,j,k,l\} = \{1,\dots,8\}\Bigg\},\] and $8$ rulings of
  the form $8H-3(E_1+\dots+E_8)+2E_i$:
    \[\{F_{i,j}+C_{j,i} \mid j \ne i\}.\]
  \item $8 \cdot \binom 7 2 = 168$ rulings of the form
    $9H-3(E_1+\dots+E_8)+E_i-(E_j+E_k)$:
    \[\{S_j+C_{k,i},S_k+C_{j,i},F_{i,l}+V_{j,k,l} \mid l \notin \{i,j,k\}\}.\]
  \item $\binom 8 4 = 70$ rulings of the form
    $10H-4(E_1+\dots+E_8)+E_i+E_j+E_k+E_l$:
    \[\{F_{i,j}+F_{k,l},F_{i,k}+F_{j,l},F_{i,l}+F_{j,k},S_a+V_{b,c,d}
    \mid \{a,b,c,d,i,j,k,l\} = \{1,\dots,8\}\}.\]
  \item $8$ rulings of the form $11H-4(E_1+\dots+E_8)+E_i$:
    \[\{S_j+F_{i,j} \mid j \ne i\}.\]
  \end{itemize}
  There is no way to write a ruling as the sum of $\anti 8$ and negative
  curves.
\end{lemma}

\begin{proof}
  Because of the Weyl group symmetry, we need to prove the last
  statement only in one case, say $H-E_1$. In this case, it is obvious.
  
  By Table~\ref{tab:relations}, there can be no other rulings, and
  each ruling can be expressed in no further ways as the sum of two
  negative curves.
\end{proof}

With this information, Lemma~\ref{lem:quadric_relations} allows us to
determine the 10800 relations coming from rulings explicitly.

We can find more quadratic relations in the degrees corresponding to
$(2)$-rulings: Because of the Weyl group symmetry, it is enough to consider
the $(2)$-ruling $D := E_2+C_{2,1}$. This can also be written as $E_j+C_{j,1}$
for any $j \ne 1$ and as $m_{i,j}+Q_{1,i,j}$ for any $i,j \ne 1$, giving 28
section in $\Gamma(S_8, D)$. As $D=\anti 8+E_1$, we get two further section
$\eta_1\kappa_1, \eta_1\kappa_2$. As the previous quadratic relations do not
induce relations in this degree of $\Cox(S_8)$, and because we calculate
$\dim\Gamma(S_8, D) = 3$ for this nef degree, we obtain 27 relations, which
can be calculated explicitly as before.

Every negative curve has intersection number $2$ with
exactly 56 other curves (e.g. $(E_1,D)=2$ if and only if $D \in
\{C_{1,i}, V_{1,i,j}, F_{i,j}, T_i\}$ for $i,j\ne 1$), so it occurs in
exactly 56 $(2)$-rulings. On the other hand, as every $(2)$-ruling can be
written in 28 ways as the sum of two negative curves, the total number
of $(2)$-rulings is $\frac{240\cdot 56}{2\cdot 28} = 240$. Therefore, we
obtain another 6480 relations from the $(2)$-rulings. To determine them
explicitly, we need the following more detailed information:

\begin{lemma}
  Each of the 240 $(2)$-rulings can be written as a sum of two negative
  curves in the following 28 ways:
  \begin{itemize}
  \item $8$ $(2)$-rulings of the form \[\anti 8+E_i=3H-(E_1+\dots+E_8)+E_i:\]
    \[\{E_j+C_{j,i}, m_{j,k}+Q_{i,j,k} \mid j,k \ne i\}.\]
  \item $\binom 8 2 = 28$ $(2)$-rulings of the form
    \[\anti 8+m_{i,j}=4H-(E_1+\dots+E_8)-(E_i+E_j):\]
    \[\Bigg\{
    \begin{aligned}
      &E_k+V_{i,j,k}, m_{i,k}+C_{j,k},\\
      &m_{j,k}+C_{i,k}, Q_{a,b,c}+Q_{d,e,f}
    \end{aligned}
    \Bigg| 
    \begin{aligned}
      &k \notin \{i,j\}, \\
      &\{i,j,a,b,c,d,e,f\} = \{1,\dots,8\}
    \end{aligned}
    \Bigg\}.\]
  \item $\binom 8 3 = 56$ $(2)$-rulings of the form
    \[\anti 8+Q_{i,j,k}=5H-2(E_1+\dots+E_8)+E_i+E_j+E_k:\]
    \[
    \Bigg\{
    \begin{aligned}
      &E_i+F_{j,k},E_j+F_{i,k},E_k+F_{i,j},m_{a,b}+V_{c,d,e},\\
      &Q_{i,j,l}+C_{l,k},Q_{i,k,l}+C_{l,j},Q_{j,k,l}+C_{l,i}
    \end{aligned}
    \Bigg| 
    \begin{aligned}
      &\{i,j,k,a,b,c,d,e\} \\ &= \{1,\dots,8\},l \notin \{i,j,k\}
    \end{aligned}
    \Bigg\}.
    \]
  \item $8\cdot 7 = 56$ $(2)$-rulings of the form
    \[\anti 8+C_{i,j}=6H-2(E_1+\dots+E_8)-E_i+E_j:\]
    \[\{E_j+T_i, m_{i,k}+F_{j,k}, Q_{j,k,l}+V_{i,k,l}, C_{i,k}+C_{k,j} 
    \mid k,l \notin\{i,j\}\}.\]
  \item $\binom 8 3 = 56$ $(2)$-rulings of the form
    \[\anti 8+V_{i,j,k}=7H-2(E_1+\dots+E_8)-(E_i+E_j+E_k):\]
    \[
    \Bigg\{
    \begin{aligned}
      &T_i+m_{j,k},S_j+m_{i,k},S_k+m_{i,j},F_{a,b}+Q_{c,d,e},\\
      &V_{i,j,l}+C_{k,l},V_{i,k,l}+C_{j,l},V_{j,k,l}+C_{i,l}
    \end{aligned}
    \Bigg|
    \begin{aligned}
      \{i,j,k,a,b,c,d,e\}\\ = \{1,\dots,8\}, l \notin \{i,j,k\}
    \end{aligned}
    \Bigg\}.
    \]
  \item $\binom 8 2 = 28$ $(2)$-rulings of the form
    \[\anti 8+F_{i,j}=8H-3(E_1+\dots+E_8)+E_i+E_j:\]
    \[\Bigg\{
    \begin{aligned}
      &S_k+Q_{i,j,k}, F_{i,k}+C_{k,j},\\&F_{j,k}+C_{k,i},
      V_{a,b,c}+V_{d,e,f}
    \end{aligned}
    \Bigg| 
    \begin{aligned}
      &k \notin \{i,j\},\\&\{i,j,a,b,c,d,e,f\} = \{1,\dots,8\}
    \end{aligned}
    \Bigg\}.\]
  \item $8$ $(2)$-rulings of the form
    \[\anti 8+T_i=9H-3(E_1+\dots+E_8)-E_i:\]
    \[\{S_j+C_{i,j}, F_{j,k}+V_{i,j,k} \mid j,k \ne i\}.\]
  \end{itemize}
\end{lemma}

Furthermore, the 242 generators give the 123 quadratic monomials
\[\eta_i\sigma_i,\quad\mu_{i,j}\psi_{i,j},\quad \nu_{i,j,k}\phi_{i,j,k},\quad
\lambda_{i,j}\lambda_{j,i},\quad \kappa_1^2, \kappa_1\kappa_2,
\kappa_2^2\] in the 4-dimensional subspace $\Gamma(S_8, \tanti 8)$ of
$\Cox(S_8)$. Note that $\tanti 8$ is the unique $(3)$-ruling. As the
relations coming from rulings and $(2)$-rulings do not induce
relations in $\Gamma(S_8, \tanti 8)$, we obtain another 119 relations.
Their equations can be calculated in the same way as before.

\begin{lemma}
  There are exactly $17399$ independent quadratic relations in $\Cox(S_8)$.
\end{lemma}

\begin{proof}
  The relations in $\Cox(S_8)$ are generated by relations which are
  homogeneous with respect to the $\Pic(S_8)$-grading. A quadratic
  relation involving a term $\xi(D_1)\xi(D_2)$ has degree $D=D_1+D_2$.
  The relations of degree $D_1+D_2$ depend on the intersection number
  $n=(D_1,D_2)$:
  \begin{itemize}
  \item If $n=1$, then $D$ is a $(1)$-ruling. As described above, we have
    exactly 10800 corresponding relations.
  \item If $n=2$, then $D$ is a $(2)$-ruling. We have described the 6480
    resulting relations.
  \item If $n=3$, then $D = \tanti 8$, which results in exactly 119
    quadratic relations.
  \item If $n=0$, then $D=D_1+D_2$ is not nef since $(D,D_1)=-1$.
    However, by results of \cite[Section 3]{MR2029868}, the relations
    in $\Cox(S_8)$ are generated by relations in nef degrees.
  \item If $n=-1$, then $D_1=D_2$, and $(D,D_1) = -2$, so $D$ is not
    nef, giving no generating relations as before.
  \end{itemize}
  There are no other quadratic relations involving $\kappa_i$ because
  the 240 degrees $\anti 8+D_1$ for some negative curve $D_1$ are
  exactly the $(2)$-rulings, and the degree $\tanti 8$ has also been
  considered.
\end{proof}

Let $J_8$ be the ideal generated by these $17399$ quadratic relations in
\[R_8 = \KK[\eta_i,\, \mu_{i,j},\, \nu_{i,j,k},\,\lambda_{i,j},\,
\phi_{i,j,k},\,\psi_{i,j},\,\sigma_i,\,\kappa_i].\]

\

\textbf{Proof of Lemma~\ref{lem:dependence}.} Let $D \in \DD_8 \setminus
\{K_1,K_2\}$ be any negative curve. We call a variable $\xi(D')$ for a
negative curve $D' \in \DD_8$ an \emph{$(n)$-variable} if $(D,D') = n$.

As for $r=7$ in the previous section, we show that the rulings
determine the $(1)$-variables in terms of the $(0)$-variables and
$\xi(D) \ne 0$.

For the two variables $\kappa_i = \xi(K_i)$ corresponding to $\anti
8$, we use the $(2)$-ruling $\anti 8+D$: As $(D, \anti 8+D) = 0$, we
have $(D, D_i)=(D,D_i')=0$ for any other possibility to write $\anti
8+D$ as the sum of two negative curves $D_i, D_i'$. Since $(\anti 8
+D, \anti 8) = 2$, by \cite[Prop.  3.4]{MR2029863}, the quadratic
monomials $\xi(D_i)\xi(D_i')$ span $\Gamma(S_8, \anti 8+D)$, so this
$(2)$-ruling induces relations of the form
\[\kappa_i \xi(D) = \sum a_i \xi(D_i)\xi(D_i').\] Therefore,
$\kappa_1, \kappa_2$ are determined by $\xi(D)$ and the
$(0)$-variables.

Any $(2)$-coordinate $\xi(D')$ is determined by the $(2)$-ruling
$D+D'$: As $(D,D+D') = 1$, we have $(D,D_i)=0$ and $(D,D_i')=1$ for
every other possibility to write $D+D'$ as the sum of two negative
curves $D_i, D_i'$. Furthermore, if $D+D' = \anti 8 + D''$, then
$(D,D'')=0$. Therefore, the relations corresponding to this
$(2)$-ruling determine $\xi(D')$ in terms of the $(0)$- and
$(1)$-variables and $\kappa_1, \kappa_2, \xi(D)$.

Finally, there is a unique $(3)$-coordinate $D'$, where $D+D' = \tanti
2$ is the $(3)$-ruling. As all other variables are known at this
point, the relations corresponding to $\tanti 8$ containing the term
$\xi(D)\xi(D')$ determine $\xi(D')$.

Consider a point in $U_{K_j}$, i.e., with $\kappa_j \ne 0$. As above, by
\cite[Prop.  3.4]{MR2029863}, $\Gamma(S_8, \tanti 8)$ is spanned by the
monomials $\xi(D_i)\xi(D_i')$ for $(3)$-rulings $D_i,D_i'$. Therefore, we have
relations of the form
\[\kappa_j^2 = \sum a_i \xi(D_i)\xi(D_i'),\] which shows that $\xi(D_i) \ne 0$
for some $i$. This proves that $Z_8 \setminus \{0\}$ is covered by the sets
$U_D$ for $D \in \DD_8 \setminus \{K_1,K_2\}$.

\

\textbf{Proof of Lemma~\ref{lem:smooth_point}.} Let $p \in \Spec
R_8$ be the point whose coordinates are zero, except that $\eta_8,
\mu_{1,3}, \mu_{2,3}, \mu_{3,4}, \mu_{3,5}, \mu_{3,6}, \mu_{3,8}$ are
$1$ and
\[(\eta_1, \eta_2, \eta_4, \eta_5, \eta_6) = 
\left(\frac{\alpha_3\alpha_4}\alpha, \frac{\alpha_4}\alpha,
\frac{(1-\alpha_3)\alpha_4}\alpha,
\frac{(\alpha_1-\alpha_3)\alpha_4}{\alpha_1\alpha},
\frac{(\alpha_2-\alpha_3)\alpha_4}{\alpha_2\alpha}\right),\] where
$\alpha:=\alpha_4-\alpha_3$.  This point satisfies the five equations
corresponding to the ruling $H-E_3$:
\begin{align*}
  &\eta_1\mu_{1,3} - \frac{\alpha_3\alpha_4}{\alpha_3 - \alpha_4}
  \eta_7\mu_{3,7} + \frac{\alpha_3\alpha_4}{\alpha_3 - \alpha_4}
  \eta_8\mu_{3,8}=0,\\
  &\eta_2\mu_{2,3} - \frac{\alpha_3}{\alpha_3 - \alpha_4}\eta_7\mu_{3,7} +
  \frac{\alpha_4}{\alpha_3 - \alpha_4}\eta_8\mu_{3,8}=0,\\
  &\eta_4\mu_{3,4} + \frac{\alpha_3\alpha_4 - \alpha_3}{\alpha_3 - \alpha_4}
  \eta_7\mu_{3,7} + \frac{-\alpha_3\alpha_4 + \alpha_4}{\alpha_3 - \alpha_4}
  \eta_8\mu_{3,8}=0,\\
  &\eta_5\mu_{3,5} + \frac{-\alpha_1\alpha_3 + \alpha_3\alpha_4}
  {\alpha_1\alpha_3 - \alpha_1\alpha_4}\eta_7\mu_{3,7} +
  \frac{\alpha_1\alpha_4 - \alpha_3\alpha_4}
  {\alpha_1\alpha_3 - \alpha_1\alpha_4}\eta_8\mu_{3,8}=0,\\
  &\eta_6\mu_{3,6} + \frac{-\alpha_2\alpha_3 + \alpha_3\alpha_4}
  {\alpha_2\alpha_3 - \alpha_2\alpha_4}\eta_7\mu_{3,7} +
  \frac{\alpha_2\alpha_4 - \alpha_3\alpha_4}{\alpha_2\alpha_3 -
    \alpha_2\alpha_4}\eta_8\mu_{3,8}=0.
\end{align*}
Consider intersection numbers between the negative curves
corresponding to the twelve non-zero coordinates. They are zero except
for the six pairs corresponding to the ruling $H-E_3$. Therefore, no
pair of non-zero coordinates occurs in relations corresponding to
other $(n)$-rulings, which shows that $p \in \AA(S_8)$. We check
directly that the Jacobian in $p$ has full rank 231.

\bibliographystyle{alpha}

\bibliography{cox_smooth}

\end{document}